\documentclass[12pt,a4paper]{article}
\usepackage{geometry}
\usepackage{fancyhdr}
\usepackage[ansinew]{inputenc} 
\usepackage[T1]{fontenc}
\usepackage{a4wide}
\usepackage[english]{babel}
\usepackage{moreverb}
\usepackage{amsfonts}
\usepackage{mathrsfs}
\usepackage{amsmath}
\usepackage{amssymb}
\usepackage{amsthm}
\usepackage{enumerate}
\usepackage[pdftex]{graphicx}
\usepackage{hyperref}
\usepackage{fancyhdr}
\usepackage{parskip}
\usepackage{color}
\usepackage[dvipsnames]{xcolor}

\usepackage{sidecap}
\usepackage{wrapfig}

\usepackage{marvosym}

\usepackage{tkz-berge}
\usepackage{tkz-graph}

\usepackage{tikz}
\usepackage{calc}
\usetikzlibrary{decorations.markings}
\tikzstyle{vertex}=[circle, draw, inner sep=0pt, minimum size=6pt]

\usepackage[makeroom]{cancel}

\renewcommand{\l}{\ell}

\newcommand{\Z}{\mathbb{Z}}

\newcommand{\R}{\mathbb{R}}

\newcommand{\N}{\mathbb{N}}

\newcommand{\U}{\mathscr{U}}

\renewcommand{\phi}{\varphi}
\renewcommand{\epsilon}{\varepsilon}

\renewcommand{\Z}{\mathbb{Z}}
\renewcommand{\P}{\mathcal{P}}

\renewcommand{\R}{\mathbb{R}}

\renewcommand{\N}{\mathbb{N}}

\renewcommand{\U}{\mathcal{U}}

\renewcommand{\bar}{\overline}

\theoremstyle{definition}\newtheorem{theorem}{Theorem}[section]
\theoremstyle{definition}\newtheorem{example}[theorem]{Example}
\theoremstyle{definition}\newtheorem{lemma}[theorem]{Lemma}
\theoremstyle{definition}\newtheorem{definition}[theorem]{Definition}
\theoremstyle{definition}\newtheorem{proposition}[theorem]{Proposition}
\theoremstyle{definition}
\theoremstyle{definition}\newtheorem{remark}[theorem]{Remark}
\theoremstyle{definition}
\theoremstyle{conjecture}

\usepackage{subcaption}

\begin{document}

\title{More on the Generalised Shift Graph}

\author{Milette Riis}

\LARGE
\center{\underline{{\textbf{The Generalised Shift Graph}}}

\normalsize

M. Riis (University of Leeds, UK)

\flushleft

\bigbreak

\section*{Abstract}

In 1968, Erd\"os\footnote{\cite{shiftnew} attributes this to \cite{erdos}} defined the Shift Graph as the graph whose vertices are the $k$-element subsets of $[n]=\{0,1,2,...,n-1\}$ such that $A=\{a_1,...,a_k\}$ and $B=\{b_1,...,b_k\}$ are neighbours iff $a_1<b_1=a_2<b_2=a_3<... <b_{n-1}=a_n<b_n$. In the paper \textit{On the Generalised Shift Graph}, Avart, Luczac and R\"odl extend this definition to include all possible arrangements of the $a_is$ and $b_is$, known as \textit{types}. In this paper, we will consider a selection of these types and study the corresponding graphs.

We are interested in to what extent the graphs $G(S,\tau)$ and $G(S',\tau)$ are distinct for distinct linear orderings $S,S'$ and for some type $\tau$. In this paper, we will concentrate on ordinals and types of the form $\sigma_{a,b}=11...133...322...2$. We will show that  if $G(\alpha,\sigma_{a,b})\cong G(\beta,\sigma_{a,b})$ then  $\alpha=\beta$. We will also consider the chromatic number and the automorphism groups of these graphs in order to gain a deeper understanding of their properties.

\section{Background Material}

We would like to describe a new class of graphs. Each graph is determined by two things: a totally ordered set $S$, and a ``type'' $\tau$; we thus denote these graphs by  $G(S,\tau)$. Both the edges and the vertices of $G$ are encoded by $\tau$.

\bigbreak

\definition Let $k,\l\in \N$ be fixed, $k\leq \l$. We say that a sequence $\tau=(\tau_i)_{i=1}^\l$ is a type of width $k$ and length $\l$ if $\tau_i\in \{1,2,3\}$ and $|\{i:\tau_i\in \{1,3\}\}|=|\{i:\tau_i\in \{2,3\}\}|=k$.\footnote{This definition is taken from \cite{shiftnew}, \textit{pp.173-174}.} 

i.e. a type is a sequence of 1s, 2s, and 3s such that the number of 1s is equal to the number of 2s. \newline

We interpret this type as follows:\newline

\definition  Let $x$ and $y$ be $k$-element subsets (listed in increasing order) of some totally ordered set $(S,<)$. Let $x\cup y=\{z_1,...,z_{\l}\}$, with $z_1<z_2<...<z_\ell$. Then we say that the pair $x, y$ \textit{has type} $\tau$ (denoted by $t(x, y)=\tau$) iff: 
\begin{align*}
\tau_i=1 &\Rightarrow z_i\in x\setminus y,\\  \tau_i=2 &\Rightarrow z_i\in y\setminus x,\text{ and }\\ \tau_i=3 &\Rightarrow z_i\in x\cap y
\end{align*}\newline

Throughout this paper, we may assume that any subset $\{x_1,...,x_n\}$ of a totally ordered set $S$ is ordered by the induced ordering.

\bigbreak

\definition The graph $G(S,\tau)$ is defined to be the graph whose vertices are the $k$-element subsets of $S$, and where there is an edge between $x$ and $y$ iff $t(x,y)=\tau$.\newline

\example Consider the graph $G(\omega, 12312)$. Then $\tau=12312$ has width 3 and length 5. Now let $x=(1,5,6)$ and $y=(3,5,8)$. Then $x\cup y=\{1,3,5,6,8\}$, and thus $t(x,y)= \tau$. However, if we let $x'=(1,3,5)$ and $y'=(2,4,6)$, then $x' \cup y'=\{1,2,3,4,5,6\}$, and so $t(x',y')\neq \tau$. In fact, $t(x',y')=121212$. Thus in the graph $G(\omega,12312)$, there is an edge between the vertices $x$ and $y$ but not between $x'$ and $y'$.

\bigbreak

\example Consider the graph $G(5,1221)$ consisting of 2-element subsets of $\{0,1,2,3,4\}$. There is an edge between any two vertices $x=(x_1,x_2)$ and $y=(y_1,y_2)$ such that $x_1<y_1<y_2<x_2$. For example, there is an edge between $(1,4)$ and $(2,3)$, but not between $(1,3)$ and $(2,4)$. \newline

\example Consider the graph $G(\R, 13332)$. Then $|\{i:\tau_i\in\{1,3\}\}|=|\{i:\tau_i\in \{2,3\}\}|=4$, and so the vertices of $G(\R,13332)$ consist of elements of $\R^4$. There is an edge between two vertices $x=(x_1,x_2,x_3,x_4)$ and $y=(y_1,y_2,y_3,y_4)$ iff $x_1<y_1=x_2<y_2=x_3<y_3=x_4<y_4$. Any graph like this with type $133...32$ is an example of the \textit{shift graph} \cite{erdos}.\newline

\textbf{Notation:} In this paper, we define $\sigma_{a,b}$ to be the type $11...133...322...2$, with $a$ copies of 1 and 2, and $b$ copies of 3. Note the following cases which are of particular interest:

$\ \ \ \ \ \bullet \ \ \sigma_{a,0}=11...122...2$

$\ \ \ \ \ \bullet \ \ \sigma_{0,b}=33...3$

$\ \ \ \ \ \bullet \ \ \sigma_{a,1}=11...1322...2$

$\ \ \ \ \ \bullet \ \ \sigma_{1,b}=133...32$

\section{Intuition}

Let $S$ be a totally ordered set as before. Our general question is to see to what extent $S$ can be recovered from $G(S,\tau)$ for any type $\tau$.  We will focus mainly on the case of $S$ an ordinal, and $\tau$ restricted to a set of possible types. We will start by considering the type $\tau=\sigma_{1,1}=132$. The graph $G(S,132)$ consists of vertices of the form $(x,y)$ with $x,y\in S$ and edges between $(x,y)$ and $(y,z)$ for all $x,y,z\in S$ where $x<y<z$.

\bigbreak

\proposition $G(S,132)$ contains an isolated point iff $S$ contains a smallest and a largest element.

\proof Suppose $S$ contains a smallest point $\alpha$ and a largest point $\beta$. Then $(\alpha,\beta)$ is an isolated point. Conversely, let $(\alpha,\beta)$ be an isolated point in $G(S,\tau)$ and suppose $\alpha$ isn't the smallest element of $S$. Then there exists some $\gamma<\alpha$, and so the point $(\gamma,\alpha)$ is joined to $(\alpha,\beta)$, a contradiction since $(\alpha,\beta)$ is isolated.  Hence $\alpha$ is the smallest element of $S$, and similarly $\beta$ is the largest element of $S$.\qed

\label{isolated}

\bigbreak

We will need the following definitions from graph theory:

\bigbreak

\definition Let $G$ be any graph. A set of vertices $S\subseteq V(G)$ is a \textit{clique} if there is an edge between every pair of vertices in $S$.

\bigbreak

\definition Let $G$ be any graph. A set of vertices $S\subseteq V(G)$ is a \textit{co-clique} if $S$ contains no edges.

\bigbreak

\definition Let $\kappa$ be any cardinal.
The complete graph $K_{\kappa}$ consists of $\kappa$ vertices such that there is an edge between every pair of vertices.

\bigbreak

\definition Let $\kappa_1,\kappa_2$ be any cardinals. The \textit{complete bipartite graph} $K_{\kappa_1,\kappa_2}$ consists of a co-clique $S_1$ of $\kappa_1$ vertices and a co-clique $S_2$ of $\kappa_2$ vertices such that every vertex in $S_1$ is joined to every vertex in $S_2$.

\bigbreak

\definition Let $v$ be a vertex in a graph $G$. Then the \textit{degree} of $v$, denoted by $d(v)$, is the number of edges incident to $v$.

\bigbreak

\lemma Let $\lambda,\kappa>1 $ be any cardinals (finite or infinite) and let $\alpha$ be an ordinal. Then  $K_{\lambda,\kappa}$  as an induced subgraph of $G(\alpha,132)$ must have the form $\{(a,x),(x,b):a\in X, b\in Y\}$ for some $X,Y$ such that $|X|=\lambda$, $|Y|=\kappa$.

\center\includegraphics[scale=0.2]{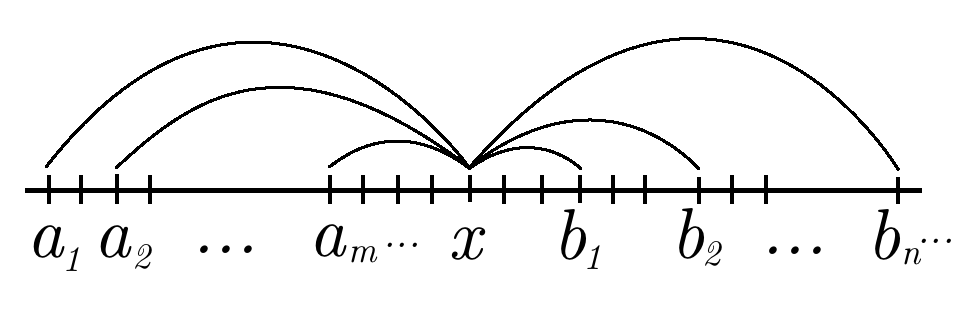}

\flushleft

\proof We will show that the result holds for $\kappa=2$ and $\kappa=3$, and thus that the result holds for all $\kappa\geq 2$. 

Case (1): $\kappa=2$. Let $(x,y)$ be a common neighbour of $(a_1,a_2),(b_1,b_2)$ where $a_1 \leq b_1$. Then $(y=a_1 \vee x=a_2)$ and $(y=b_1\vee x=b_2)$. If $y=a_1$, then since $x<y=a_1\leq b_1<b_2$, we must have $y=b_1$ and so $a_1=b_1$.
If $x=a_2$ and $y=b_1$ then $a_1<a_2<b_1<b_2$, and clearly $(a_1,a_2)$ and $(b_1,b_2)$ have exactly one common neighbour. Thus if $x=a_2$ and $y\neq b_1$, we must have $x=b_2$ i.e. $a_2=b_2$. Hence for 2  vertices $(a_1,a_2),(b_1,b_2)$ to share $\lambda>1$ many neighbours, we must have $a_1=b_1$ or $a_2=b_2$.

Case (2): $\kappa=3$. For a set of 3 vertices $(a_1,a_2),(b_1,b_2),(c_1,c_2)$ to share a set of $\lambda$ neighbours, we know that for each pair we must have $x_i=y_i$ for some $x,y\in \{a,b,c\}, i\in \{1,2\}$. We would like to show that we must have either $a_1=b_1=c_1$ or $a_2=b_2=c_2$. Suppose without loss of generality that $a_1=b_1$. If $c_1\neq a_1$, then we must have $c_2=b_2$ and $c_2=a_2$. Similarly, suppose  $a_2=b_2$. If $c_2\neq a_2$, then we must have $c_1=b_1$ and $c_1=a_1$. Hence either $a_1=b_1=c_1$ or $a_2=b_2=c_2$.

Now let $\lambda>1$ be any cardinal, and suppose we have a set $S$ of $\kappa>1$ vertices. Look at the first two vertices $(a_1,a_2)$ and $(b_1,b_2)$ of $S$. By Case (1), either $a_1=b_1$ or $a_2=b_2$. Without loss of generality assume $a_1=b_1$. Now let $(x_1,x_2)$ be any third vertex in $S$. By Case (2), we must have either  $a_1=b_1=x_1$ or $a_2=b_2=x_2$. Assuming the vertices are distinct, our assumption implies that  $x_1=a_1=b_1$. Thus all the first coordinates of the vertices in $S$ are equal. A similar argument applies if $a_2=b_2$ for the vertices $(a_1,a_2)$ and $(b_1,b_2)$.
\qed

\label{completebipartite}

\bigbreak

\proposition $G(\omega+x,132)\not\cong G(\omega+y,132)$ for finite $x\neq y$.

\proof Using a simple counting argument, we see that the vertex $(n,\omega+x-m-1)$ has degree ${n}+{m}$ for all $n\geq 0$ and $m<x$. Thus the number of vertices of degree $d$ equals the number of pairs $\langle n,m \rangle$ such that ${n}+{m}=d$.

Hence, by looking at the degrees of the vertices of the graph $G(\omega+x,132)$, we can determine  the value of $x$.\qed

\label{2.8}

We can similarly show that $G(\omega+x,\sigma_{a,1})\not\cong G(\omega+y,\sigma_{a,1})$ for finite $x\neq y$ and $a>0$. Notice that if $\omega$ is replaced by any limit ordinal $\alpha$, Proposition \ref{2.8} remains valid.

\pagebreak

\proposition  $G(\omega,132)\ncong G(\omega+\omega,132)$.

\proof We will show that $G(\omega+\omega,132)$ contains a copy of $K_{\omega,\omega}$, whereas $G(\omega,132)$ does not.

In $G(\omega+\omega,132)$, the set of vertices $\{(x,\omega),(\omega,\omega+y):x,y<\omega, y\neq 0\}$ forms a copy of $K_{\omega,\omega}$.

\textit{Claim:} Such a set cannot exist in $G(\omega,\tau)$.

\textit{Proof of Claim:} By Lemma \ref{completebipartite}, a set of $\omega$ many points which all share at least two neighbours must also either begin with the same ordinal or end with the same ordinal. In the graph $G(\omega,132)$, it is impossible for $\omega$ many points to all end with the same ordinal. Thus a set of $\omega$ many points which all share a neighbour in $G(\omega,\tau)$ must be of the form $F_x=\{(x,y):y>x\}$ for some fixed $x$. The shared neighbour set of $F_x$ must therefore be finite as  $x$ is finite. Hence $G(\omega,\tau)$ does not contain a copy of $K_{\omega,\omega}$.
\flushright
 \textit{Q.E.D. \ \ Claim} 
\flushleft \qed

\proposition  $G(\omega+\omega,132)\ncong G(\omega+\omega+\omega,132)$.

\textit{Sketch of Proof.} We will show that $G(\omega+\omega+\omega,132)$ contains two copies of $K_{\omega,\omega}$ joined together by pairs of vertices by a matching, whereas $G(\omega+\omega,132)$ does not.

\center

\includegraphics[scale=0.12]{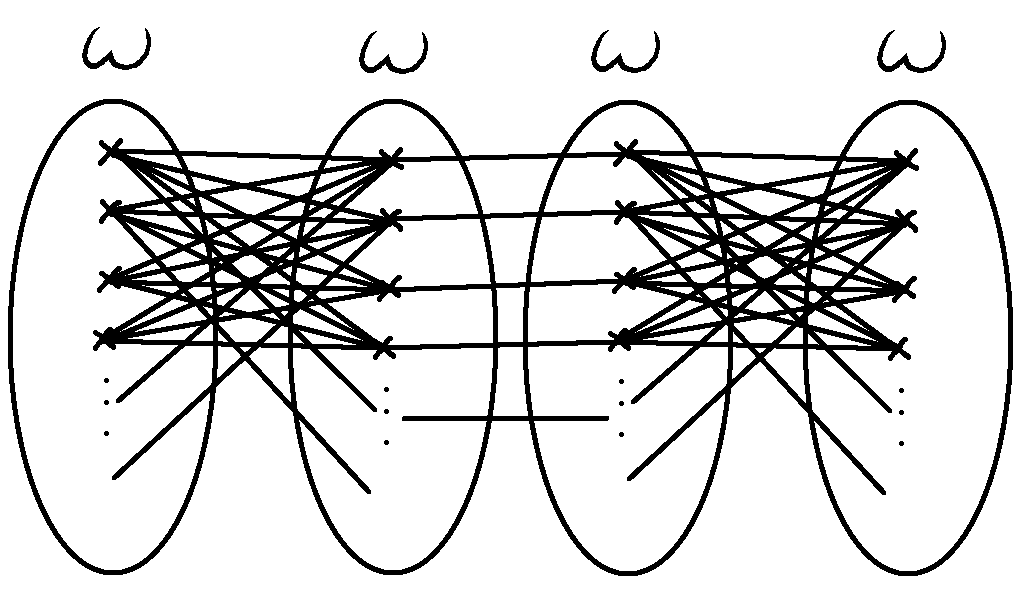}

\flushleft

In $G(\omega+\omega+\omega,132)$, the set of vertices $\{(x,\omega),(\omega,\omega+y):x,y<\omega\}$ forms a copy of $K_{\omega,\omega}$, along with the set of vertices  $\{(\omega+x,\omega+\omega) , (\omega+\omega,\omega+\omega+y):x,y<\omega\}$. Every vertex in $\{(\omega,\omega+y):y<\omega\}$ is joined to exactly one vertex in $\{(\omega+x,\omega+\omega):x<\omega\}$, forming the desired subgraph.

Such a subgraph cannot exist in $G(\omega+\omega,132)$, since the copies of $K_{\omega,\omega}$ must have the form $\{(x,\omega+n),(\omega+n,\omega+y):x<\omega,n<y<\omega\}$ and $\{(x,\omega+m),(\omega+m,\omega+y):x<\omega,m<y<\omega\}$ for some fixed $m,n$. Without loss of generality, assume $m>n$. Then the vertex $(\omega+n,\omega+m)$ lies in both copies of $K_{\omega,\omega}$. \qed

\bigbreak

It follows from the proof of this theorem that $G(\omega\cdot n,132)\ncong G(\omega \cdot n',132)$ for $n<n'$, as we can expand the proof to include $n'-1$ copies of $K_{\omega,\omega}$ joined as above by a matching.\newline
It also follows that $G(\omega_k \cdot n,132)\ncong G(\omega_k \cdot n',132)$, as we can repeat the same argument with different ordinals $\omega_k$. Thus $G(\omega_k \cdot n,132)\ncong G(\omega_{\ell} \cdot n',132)$ in general for $k\neq \ell,\ n<n'$.

\bigbreak

\bigbreak

We will now consider the type $\tau=\sigma_{2,0}=1122$. If $\kappa$ is any cardinal, then the graph $G(\kappa,1122)$ consists of all vertices of form $(x,y)$ with $x,y\in \kappa$ and with an edge between $(x,y)$ and $(z,w)$ (with $x<z$) iff $y<z$ for all $x,y,z,w\in \kappa$.

Cut $\kappa$ into $\kappa$ disjoint sets $A_{\alpha}$ for $\alpha<\kappa$. Consider a set of pairwise disjoint pairs $(a_{\beta},b_{\beta})$ in $A_{\alpha}$, $a_0<b_0<a_1<b_1<...$. This gives rise to a copy of $K_{\kappa}$ inside each $A_{\alpha}$, thus $\kappa$ disjoint copies of $K_{\kappa}$.

\bigbreak

\proposition $G(\omega+x,1122)\not\cong G(\omega+y,1122)$ for finite $x\neq y$.

\proof Using a simple counting argument, we see that the vertex $(n,\omega+x-m-1)$ has degree ${n\choose 2}+{m\choose 2}$ for all $n\geq 0$ and $m<x$. Thus the number of vertices of degree $d$ equals the number of pairs $\langle n,m \rangle$ such that ${n\choose 2}+{m\choose 2}=d$.

Hence, by looking at the degrees of the vertices of the graph $G(\omega+x,1122)$, we can determine  the value of $x$.\qed

\label{2.11}

\bigbreak

We can similarly show that  $G(\omega+x,\sigma_{a,0})\not\cong G(\omega+y,\sigma_{a,0})$ for finite $x\neq y$ and $a>1$. Notice that if $\omega$ is replaced by any limit ordinal $\alpha$, Proposition \ref{2.11} remains valid.

\bigbreak

\proposition $G(\omega,1122)\not\cong G(\omega+\omega,1122)$.

\proof We will show that $G(\omega+\omega,1122)$ contains a copy of $K_{\omega,\omega}$ while $G(\omega,1122)$ does not. 

First we see that the set of vertices in $G(\omega+\omega,1122)$ of the form $\{(0,x):\omega>x>0\}\cup \{(\omega,\omega+x):\omega>x>0\}$ give a copy of $K_{\omega,\omega}$. We would now like to show that such a set cannot exist in $G(\omega,1122)$. 

\textit{Claim:} Any infinite co-clique in $G(\omega,1122)$ must contain an infinite subset of $\bar{x}=\{(x,y):y>x\}$ for some $x$.

\textit{Proof of Claim:} Let $(a,b)$ lie in an infinite co-clique $S$ in $G(\omega,1122)$, and let $(x_i,y_i)$ be another member of the co-clique for all $i\in \omega$. Since $S$ is a co-clique, we must have $x_i\leq b$ for infinitely many $i$. Since $b$ is finite, by the pigeonhole principle $x_i$ is constant on infinitely many values of $i$, and so $S$ contains an infinite subset of  $\bar{x_i}=\{(x_i,y_i):y_i>x_i\}$

\flushright
 \textit{Q.E.D. \ \ Claim} 
\flushleft

Now, any two infinite co-cliques $X,Y\in G(\omega,1122)$ with some subset of $\bar{x}\in X$ and some subset of $\bar{y}\in Y$ can at most have $y-x$ edges between them, which is finite as $x$ and $y$ are both finite. Thus we cannot obtain the complete bipartite graph $K_{\omega,\omega}$ as a subgraph of $G(\omega,1122)$.\qed

\bigbreak

\section{Stronger Results}

Ideally, we would like a stronger result than this. We will now consider the more general case, where $G=G(\alpha,\sigma_{a,b})$ for any ordinal $\alpha$.

\bigbreak

\theorem Let  $\alpha$ and $\beta$ be ordinals. Then \begin{align*}
G(\alpha,132)\cong G(\beta,132) \Rightarrow \alpha = \beta
\end{align*}
Furthermore, if $\alpha=\beta$ the isomorphism between $G(\alpha,132)$ and $G(\beta,132)$ is unique.  

\label{132iso}

\proof If $\alpha, \beta$ are finite, the result is trivial as if $\alpha\neq \beta$ the graphs have different numbers of vertices. To handle the infinite case of the first part of the theorem, it is sufficient to show that $\alpha$ can be determined from the graph $G(\alpha,132)$.

\bigbreak

Assume $\alpha=\alpha_0+k$, where $\alpha_0\neq 0$ is a limit ordinal and $k$ is finite. We start by determining $k$; this can be done by looking at the degrees of the vertices of $G$.

\textit{Claim 0:} $k$ is the largest finite number such that for some finite $j$ there are exactly $k$ vertices of degree $j$. If no such number exists, then $k=0$.

\textit{Proof of Claim 0:} For finite $n$, the vertex $(n,\alpha_0+k-1)$ has degree $n$ for all $n\geq 0$. If $n$ were infinite, the degree would also be infinite. In general, the vertex $(n,\alpha_0+k-m-1)$ has degree $n+m$ for all $n< \omega$. Counting the total number of vertices of each degree, we see that there are
\begin{align*}
j+1 \text{ vertices of degree } &j \text{ for all } j< k-1\\
k \text{ vertices of degree } &j \text{ for all } j\geq k-1\\
\end{align*}
\flushright
 \textit{Q.E.D. \ \ Claim 0} 
\flushleft

\bigbreak

We would now like to find $\alpha_0$. We start by removing every vertex of finite degree from $G=G(\alpha,132)=G(\alpha_0+k,132)$, leaving us with a new graph $G_0$. Note that each vertex of finite degree must have form $(n,\alpha_0+m)$ for some $n,m<\omega$.

Now consider a set of vertices $\bar{0}\subseteq V(G_0)$ with the following properties:
\begin{enumerate}[(i)]
\item The induced subgraph on $\bar{0}$ is a maximal co-clique
\item All neighbour sets of the inducted subgraph on $\bar{0}$ are pairwise disjoint 
\end{enumerate}

\textit{Claim 1:}  $\{(0,\beta):\beta>0\}$ is the unique set $\bar{0}$ fulfilling these two conditions.

\textit{Proof of Claim 1:} First notice that $\{(0,\beta):\beta>0\}$ satisfies (i) and (ii). Conversely, suppose $\bar{0}$ is a set fulfilling (i) and (ii), and $\bar{0}\not\subseteq \{(0,\beta):\beta>0\}$. Then $\bar{0}$ contains some vertex $(\delta_1,\delta_2)$ with $\delta_1\neq 0$. Since $\bar{0}$ is a maximal co-clique, it must also contain the vertex $(0,\delta_2)$; but then $(\delta_1,\delta_2)$ and $(0,\delta_2)$ are  two points in $\bar{0}$ whose neighbour sets are not pairwise disjoint (namely they both contain the vertex $(\delta_2,\delta_3)$ for some $\delta_3>\delta_2$), a contradiction. Note that since we are working in $G_0$, $\delta_2<\alpha_0$ and so such a $\delta_3$ always exists. If $\bar{0}\subset \{(0,\beta):\beta>0\}$ then $\bar{0}$ would not be maximal, and so $\bar{0}= \{(0,\beta):\beta>0\}$.
\flushright
 \textit{Q.E.D. \ \ Claim 1} 
\flushleft

We have defined $G_0$  as $G\backslash \{(x,y):d((x,y))<\omega\}$. We now define $G_{\lambda}$ recursively as follows for every limit $\lambda$: remove each vertex of finite degree from $G\backslash \bigcup_{\zeta<\lambda}G_{\zeta}$, leaving us with a new graph $G_{\lambda}$. Note that we do not define $G_\gamma$ where $\gamma$ is a successor ordinal, as if $\gamma=\lambda+k'$, then $G_{\gamma}=G_{\lambda}$. Thus at stage $\lambda$, we
have removed all vertices of the form $(\delta,\alpha_0+n)$ for $0\leq n<k$ for $\delta \geq \lambda$.

\bigbreak

Let $\gamma \in \alpha_0$ be such that $\gamma=\lambda+k'$ for some limit ordinal $\lambda$ and $k'<\omega$. Now consider a set of vertices $\bar{\gamma}\subseteq V(G_\lambda)$ with the following properties:
\begin{enumerate}[(i)]
\item The induced subgraph on $\bar{\gamma}$ is a maximal co-clique in $V(G_{\lambda})\backslash \bigcup_{\xi<\gamma}\bar{\xi}$ 
\item All neighbour sets of the induced subgraph on $\bar{\gamma}$  are pairwise disjoint
\end{enumerate}

\textit{Claim 2:} $\bar{\gamma}$ is uniquely defined as $\{(\gamma,\beta):\beta>\gamma \}$.

\textit{Proof of Claim 2:}  First notice that $\{(\gamma,\beta):\beta>\gamma\}$ satisfies (i) and (ii). Suppose $\bar{\gamma}\not\subseteq \{(\gamma,\beta):\beta>\gamma \}$. Then $\bar{\gamma}$ contains some vertex $(\delta_1,\delta_2)$ with $\delta_1>\gamma$. Since $\bar{\gamma}$ is a maximal co-clique, it must also contain the vertex $(\gamma,\delta_2)$; but then $(\delta_1,\delta_2)$ and $(\gamma,\delta_2)$ are  two points in $\bar{\gamma}$ whose neighbour sets are not pairwise disjoint (namely they both contain the vertex $(\delta_2,\delta_3)$ for some $\delta_3>\delta_2$), a contradiction. Again, note that since we are working in $G_{\lambda}$, $\delta_2<\alpha_0$ and so such a $\delta_3$ always exists. If $\bar{\gamma}\subset \{(\gamma,\beta):\beta>\gamma \}$ then $\bar{\gamma}$ would not be maximal, and so $\bar{\gamma}= \{(\gamma,\beta):\beta>\gamma \}$.
\flushright
\flushright
 \textit{Q.E.D. \ \ Claim 2} 
\flushleft

\bigbreak

We  see that $\alpha_0$ is the least ordinal such that $V(G_{\alpha_0})= \bigcup_{\beta<\alpha_0}\bar{\beta}$.  Thus, we can determine $\alpha=\alpha_0+k$ from the graph $G$.

\bigbreak

We will now show that  if $\alpha=\beta$ the isomorphism between $G(\alpha,132)$ and $G(\beta,132)$ is unique. Let $f$ be an isomorphism from $G(\alpha,132)$ to $G(\beta,132)$. 

\textit{Claim 3:} $f$ is uniquely determined.

\textit{Proof of Claim 3:} 
Let $v_0\in G(\alpha,132)$. Since $v_0=(x_1,x_2)$ for some $x_1,x_2\in \alpha$, $v_0$ lies in $\bar{x_1}$ and  $N_{v_0} \cap [ V(G(\alpha,132))\backslash \bigcup_{\xi<x_1} \bar{\xi}  ] =\bar{x_2}$. Now, $f(v_0)$  lies in some $\bar{\delta_1}$, and  $N_{f(v_0)}\cap [V(G(\beta,132))\backslash \bigcup_{\xi<\delta_1} \bar{\xi}] =\bar{\delta_2}$ in  for some $\delta_2$. Since $\bar{x}$ is uniquely determined for each $x<\alpha$, it follows that $x_1=\delta_1$ and $x_2=\delta_2$. Hence $f$ is uniquely determined.
\flushright
 \textit{Q.E.D. \ \ Claim 3}  
\flushleft

\qed 

What \textit{Claim 3} is essentially saying is that given the  graph $G=G(\alpha,132)$, for each vertex $v\in G$ we can determine which unique pair $(x_1,x_2)\in \alpha^2$ ``generated'' $v$.

\bigbreak

We can expand this theorem to the following with very little alteration:

\theorem Let $\tau= \sigma_{1,b}=133...32$, and let $\alpha$ and $\beta$ be ordinals. Then \begin{align*} 
G(\alpha,\sigma_{1,b})\cong G(\beta,\sigma_{1,b}) \Rightarrow \alpha = \beta
\end{align*}
Furthermore, if $\alpha=\beta$ the isomorphism between $G(\alpha,\sigma_{1,b})$ and $G(\beta,\sigma_{1,b})$ is unique.  

\proof If $\alpha,\beta$ are finite, the result is trivial as if $\alpha\neq \beta$ the graphs have different numbers of vertices. To handle the infinite case of the first part of the theorem, it is sufficient to show that $\alpha$ can be determined from the graph $G(\alpha,\sigma_{1,b})$.

\bigbreak

Assume $\alpha=\alpha_0+k$, where $\alpha_0\neq 0$ is a limit and $k$ is finite. We start by determining $k$ as in Theorem \ref{132iso}. We now construct $G_\lambda$ recursively as in Theorem \ref{132iso}, and consider a set of vertices $\bar{0}\subseteq V(G_0)$ with the following properties:
\begin{enumerate}[(i)]
\item The induced subgraph on $\bar{0}$ is a maximal co-clique
\item All neighbour sets of the induced subgraph on $\bar{0}$ are pairwise disjoint 
\end{enumerate}

\textit{Claim 1:} $\bar{0}$ is uniquely defined as $\{(0,\beta_1,\beta_2,...,\beta_b):\beta_i>0\}$.

The proof of Claim 1 is similar to the proof of \textit{Claim 1} in Theorem \ref{132iso}.

\bigbreak

Let $\gamma \in \alpha_0$ be such that $\gamma=\lambda+k'$ with $\lambda$ a limit ordinal and $k'<\omega$. Now consider a set of vertices $\bar{\gamma}\subseteq V(G_\lambda)$ with the following properties:
\begin{enumerate}[(i)]
\item The induced subgraph on $\bar{\gamma}$ is a maximal co-clique in $V(G_{\lambda})\backslash \bigcup_{\xi<\gamma} \bar{\xi}$ 
\item All neighbour sets of the induced subgraph on $\bar{\gamma}$  are pairwise disjoint
\end{enumerate}

\textit{Claim 2:} $\bar{\gamma}$ is uniquely defined as $\{(\gamma,\beta_1,\beta_2,...,\beta_b):\beta_i>\gamma \}$.

The proof of Claim 2 is similar to the proof of \textit{Claim 2} in Theorem \ref{132iso}.

\bigbreak

Then $\alpha_0$ is the least ordinal such that $V(G_{\alpha_0})=\bigcup_{\beta<\alpha_0}\bar{\beta}$.
Thus we  can determine $\alpha=\alpha_0+k$ from the  graph $G$.

\bigbreak

We will now show that  if $\alpha=\beta$ the isomorphism between $G(\alpha,\sigma_{1,b})$ and $G(\beta,\sigma_{1,b})$ is unique.   Let $f$ be an isomorphism from $G(\alpha,\sigma_{1,b})$ to $G(\beta,\sigma_{1,b})$. 

\textit{Claim 3:} $f$ is uniquely determined.

\textit{Proof of Claim 3:} 
Let $v_0\in G(\alpha,\sigma_{1,b})$. Then $v_0=(x_1,x_2,...,x_b)$ for some $x_1,...,x_b\in \alpha$, and thus $v_0$ lies in $\bar{x_1}$. If $v_1$ is a neighbour of $v_0$ lying in $V(G(\beta,\sigma_{1,b}))\backslash \bigcup_{\xi<x_1} \bar{\xi}$ then $v_1\in \bar{x_2}$. Additionally, if $v_{i+1}$ is a neighbour of $v_i$ lying in $V(G(\beta,\sigma_{1,b}))\backslash \bigcup_{\xi<x_{i+1}} \bar{\xi}$, then $v_{i+1}\in x_{i+2}$ for all $i\leq b$.

Now, $f(v_0)$ also lies in some $\bar{\delta_1}$, and if $f(v_1)$ is a neighbour of $f(v_0)$ lying in $V(G(\beta,\sigma_{1,b}))\backslash \bigcup_{\xi<\delta_1} \bar{\xi}$ then let $f(v_1)\in \bar{\delta_2}$ for some $\delta_2$.  Note that $f(v_1)$ and $f(v_0)$ are neighbours in $G(\beta,\sigma_{1,b})$ as $v_1$ and $v_0$ are neighbours in $G(\alpha,\sigma_{1,b})$. In general, if $f(v_{i+1})$ is a neighbour of $f(v_i)$ lying in $V(G(\beta,\sigma_{1,b}))\backslash \bigcup_{\xi<\delta_{i+1}} \bar{\xi}$, then let $f(v_{i+1})\in \delta_{i+2}$ for all $i\leq c$.
 Since $\bar{\gamma}$ is uniquely determined for each $\gamma<\alpha$, it follows that $\gamma_i=\delta_i$ for all $i\leq b$. Hence $f$ is uniquely determined.
 \flushright
 \textit{Q.E.D. \ \ Claim 3}  
\flushleft

\qed

\bigbreak

\label{isomorphism}

Once again, \textit{Claim 3} tells us that given the  graph $G=G(\alpha,\sigma_{1,b})$ and $v\in G$, we can determine the unique $(b+1)$-tuple which ``generated'' $v$.

We will see that \textit{Claim 3} does not hold in general for $\sigma_{a,b}$ with $a>1$.

\bigbreak

\definition Let $v$ be some vertex in a graph $G$. Then the neighbour set  $N_v$ of $v$ is defined as the set of neighbours of $v$ in $G$. 

\bigbreak

\theorem Let  $\alpha$ and $\beta$ be  ordinals. Then \begin{align*}
G(\alpha,11322)\cong G(\beta,11322) \Rightarrow \alpha = \beta
\end{align*}

\proof If $\alpha,\beta$ are finite, the result is trivial as if $\alpha\neq \beta$ the graphs have different numbers of vertices. To handle the infinite case, it is enough to show that $\alpha$ can be determined from the unlablled graph $G(\alpha,11322)$.

Assume $\alpha=\alpha_0+k$, where $\alpha_0$ is a limit ordinal and $k$ is finite. We can now determine $k$  by looking at the degrees of the vertices, as the vertex $(n,\alpha_0+k-m-1)$ has degree ${n\choose 2}+{m\choose 2}$ for all $\omega>n\geq 0$ and $k>m\geq 0$.

\bigbreak

We would now like to find $\alpha_0$. We start by defining $G_\lambda$ recursively as in Theorem \ref{132iso}. Consider a set of vertices $\bar{01}\subseteq V(G_0)$ with the following properties:
\begin{enumerate}[(i)]
\item The induced subgraph on $\overline{01}$ is a maximal co-clique
\item For any two vertices $v_1,v_2\in\overline{01}$, either $N_{v_1}=N_{v_2}$ or $N_{v_1}\cap N_{v_2} =\emptyset$.
\end{enumerate}

\textit{Claim 1:} $\bar{01}$ is uniquely defined as $\{(0,\beta_1,\beta_2):\beta_i>0\}\cup \{(1,\beta_1,\beta_2):\beta_i>1\}$.

\textit{Proof of Claim 1:}  First notice that $\{(0,\beta_1,\beta_2):\beta_i>0\}\cup \{(1,\beta_1,\beta_2):\beta_i>1\}$ satisfies (i) and (ii). Suppose $\bar{01}\not\subseteq \{(0,\beta_1,\beta_2):\beta_i>0\}\cup \{(1,\beta_1,\beta_2):\beta_i>1\}$. Then $\bar{01}$ must contain some vertex $(\delta,\beta_1,\beta_2)$ with $\delta\neq 0,1$. Now, since $\bar{01}$ is a maximal co-clique, it must contain $(0,\beta_1,\beta_2)$; but since $\delta\geq 2$ , $N_{(\delta,\beta_1,\beta_2)}$  and $N_{(0,\beta_1,\beta_2)}$ are not equal as $(0,1,\delta)$ lies in $N_{(\delta,\beta_1,\beta_2)}$ and not in $N_{(0,\beta_1,\beta_2)}$, yet have a non-empty intersection as $(\beta_2,\beta_2+1,\beta_2+2)$ lies in both, contradicting the second condition above. Note that since we are working in $G_0$, $\beta_2<\alpha_0$ and so  $\beta_2+1,\beta_2+2 \in \alpha$. If $\bar{01}\subset \{(0,\beta_1,\beta_2):\beta_i>0\}\cup \{(1,\beta_1,\beta_2):\beta_i>1\}$ then $\bar{01}$ would not be maximal, and so $\bar{01}= \{(0,\beta_1,\beta_2):\beta_i>0\}\cup \{(1,\beta_1,\beta_2):\beta_i>1\}$.
\flushright
 \textit{Q.E.D. \ \ Claim 1} 
\flushleft

Now consider a set of vertices $\bar{2}\subseteq V(G_0)$ with the following properties:
\begin{enumerate}[(i)]
\item The induced subgraph on $\bar{2}$ is a maximal co-clique in $V(G_0)\backslash\overline{01}$
\item For any two vertices $v_1,v_2\in\bar{2}$, either $N_{v_1}=N_{v_2}$ or $N_{v_1}\cap N_{v_2} =\emptyset$.
\end{enumerate}

\textit{Claim 2:} $\bar{2}$ is uniquely defined as $\{(2,\beta_1,\beta_2):\beta_i>2\}$.

\textit{Proof of Claim 2:} First notice that $\{(2,\beta_1,\beta_2):\beta_i>2\}$ satisfies (i) and (ii). Suppose $\bar{2}\not\subseteq \{(2,\beta_1,\beta_2):\beta_i>2\}$. Then $\bar{2}$ must contain some vertex $(a,\beta_1,\beta_2)$ with $a> 2$. Now, since $\bar{2}$ is a maximal co-clique, it must contain $(2,\beta_1,\beta_2)$; but then $N_{(a,\beta_1,\beta_2)}$  and $N_{(2,\beta_1,\beta_2)}$ are not equal as $(0,1,a)$ lies in $N_{(a,\beta_1,\beta_2)}$ and not in $N_{(2,\beta_1,\beta_2)}$, yet have a non-empty intersection as $(\beta_2,\beta_2+1,\beta_2+2)$ lies in both, contradicting the second condition above. Note that since we are working in $G_0$, $\beta_2<\alpha_0$ and so  $\beta_2+1,\beta_2+2 \in \alpha$. If $\bar{2}\subset\{(2,\beta_1,\beta_2):\beta_i>2\}$ then $\bar{2}$ would not be maximal, hence $\bar{2}=\{(2,\beta_1,\beta_2):\beta_i>2\}$.
\flushright
 \textit{Q.E.D. \ \ Claim 2} 
\flushleft

Let $\gamma \in \alpha_0$ be such that $\gamma=\lambda+k'$ for some limit ordinal $\lambda$ and $k'<\omega$. Consider a set of vertices $\bar{\gamma}\subseteq V(G_\lambda)$ with the following properties:
\begin{enumerate}[(i)]
\item The induced subgraph on $\bar{\gamma}$ is a maximal co-clique in $V(G_{\lambda})\backslash \bigcup_{\gamma<\lambda} \bar{\gamma}$
\item For any two vertices $v_1,v_2\in\bar{\gamma}$, either $N_{v_1}=N_{v_2}$ or $N_{v_1}\cap N_{v_2} =\emptyset$.
\end{enumerate}

\textit{Claim 3:} $\bar{\gamma}$ is uniquely defined as $\{(\gamma,\beta_1,\beta_2):\beta_i>\gamma\}$.

\textit{Proof of Claim 3:} First notice that $\{(\gamma,\beta_1,\beta_2):\beta_i>\gamma\}$ satisfies (i) and (ii). Suppose $\bar{\gamma}\not\subseteq \{(\gamma,\beta_1,\beta_2):\beta_i>\gamma\}$. Then $\bar{\gamma}$ must contain some vertex $(\delta,\beta_1,\beta_2)$ with $\delta> \gamma$. Now, since $\bar{\gamma}$ is a maximal co-clique, it must contain $(\gamma,\beta_1,\beta_2)$; but since $\delta > \gamma$ , $N_{(\delta,\beta_1,\beta_2)}$  and $N_{(\gamma,\beta_1,\beta_2)}$ are not equal as $(0,1,\gamma)$ lies in $N_{(\gamma,\beta_1,\beta_2)}$ and not in $N_{(\delta,\beta_1,\beta_2)}$, yet have a non-empty intersection as $(\beta_2,\beta_2+1,\beta_2+2)$ lies in both, contradicting the second condition above. Note that since we are working in $G_{\lambda}$, $\beta_2<\alpha_0$ and so $\beta_2+1,\beta_2+2 \in \alpha$. If $\bar{\gamma}\subset\{(\gamma,\beta_1,\beta_2):\beta_i > \gamma\}$ then $\bar{\gamma}$ would not be maximal, and so $\bar{\gamma}=\{(\gamma,\beta_1,\beta_2):\beta_i>\gamma\}$.
\flushright
 \textit{Q.E.D. \ \ Claim 3} 
\flushleft

Then $\alpha_0$ is the least ordinal such that $V(G_{\alpha_0})= \bigcup_{\beta<\alpha_0}\bar{\beta}$. Thus we can determine $\alpha=\alpha_0+k$ from the graph $G$. \qed

\label{11322}

\bigbreak

We see that we can identify the set of vertices beginning with each ordinal. We can also identify the set of vertices ending with each ordinal, as we can consider $N_v\cap G\backslash \bar{u}:u<v$ for each $v\in G(\alpha,11322)$. However, we cannot determine which 3-tuple ``generated'' some vertex $v\in G$. For example, in $G(\omega,11322)$, we cannot distinguish  say the vertices $(5,7,10)$ and $(5,8,10)$ as they have the exact same neighbour sets.

\bigbreak

\theorem Let  $\alpha$ and $\beta$ be infinite ordinals, and let $\sigma_{a,1}=11...1322...2$. Then \begin{align*}
G(\alpha,\sigma_{a,1})\cong G(\beta,\sigma_{a,1}) \Rightarrow \alpha = \beta
\end{align*}

\proof If $\alpha,\beta$ are finite, the result is trivial as if $\alpha\neq \beta$ then the graphs have different numbers of vertices. We will show that given $G=G(\alpha,\sigma_{a,1})$, we can work out what $\alpha$ is.

Assume $\alpha=\alpha_0+k$, where $\alpha_0$ is a limit and $k$ is finite. We can now work out what $k$ is by looking at the degrees of the vertices, as $G(\alpha_0+k)$ has one vertex of degree ${n\choose a} + {m\choose a}$ for all $n\geq 0$ and $k>m\geq 0$. As there are exactly $k$ vertices of degree ${j\choose a}$ for $j\geq k-1$, we see that
$k$ is the greatest finite number such that for some finite $j$ there are exactly $k$ vertices of degree ${j\choose a}$ by a similar argument to Theorem \ref{132iso} \textit{Claim 0}. If no such number exists, then $k=0$.

\bigbreak

We would now like to find $\alpha_0$. We start by definite $G_{\lambda}$ recursively as in Theorem \ref{132iso}. Consider a set of vertices $\bar{0(a-1)}\subseteq V(G_0)$ with the following properties:
\begin{enumerate}[(i)]
\item The induced subgraph on $\overline{0(a-1)}$ is a maximal co-clique
\item For any two vertices $v_1,v_2\in\overline{0(a-1)}$, either $N_{v_1}=N_{v_2}$ or $N_{v_1}\cap N_{v_2} =\emptyset$.
\end{enumerate}

\textit{Claim 1:} $\overline{0(a-1)}$ is uniquely defined as $\{(0,\beta_1,\beta_2,...,\beta_a):\beta_i>0\}\cup \{(1,\beta_1,\beta_2,...,\beta_a):\beta_i>1\} \cup ... \cup \{(a-1,\beta_1,\beta_2,...,\beta_a):\beta_i>a-1\}$.

The proof of Claim 1 is similar to the proof of \textit{Claim 1} in Theorem \ref{11322}.

\bigbreak

Now consider a set of vertices $\bar{a}\subseteq V(G_0)$ with the following properties:
\begin{enumerate}[(i)]
\item The induced subgraph on $\bar{a}$ is a maximal co-clique in $V(G_0)\backslash\overline{0(a-1)}$
\item For any two vertices $v_1,v_2\in\overline{0(a-1)}$, either $N_{v_1}=N_{v_2}$ or $N_{v_1}\cap N_{v_2} =\emptyset$.
\end{enumerate}

\textit{Claim 2:} $\bar{a}$ is uniquely defined as $\{(a,\beta_1,\beta_2,...,\beta_a):\beta_i>a\}$.

The proof of Claim 2 is similar to the proof of \textit{Claim 2} in Theorem \ref{11322}.

\bigbreak

Let $\gamma \in \alpha_0$ be such that $\gamma=\lambda+k'$ for some limit ordinal $\lambda$ and $k'<\omega$. Consider a set of vertices $\bar{\gamma}\subseteq V(G_\lambda)$ with the following properties:
\begin{enumerate}[(i)]
\item The induced subgraph on $\bar{\lambda}$ is a maximal co-clique in $V(G_{\lambda})\backslash \bigcup_{\xi<\gamma}\bar{\xi}$
\item For any two vertices $v_1,v_2\in\overline{\gamma}$, either $N_{v_1}=N_{v_2}$ or $N_{v_1}\cap N_{v_2} =\emptyset$.
\end{enumerate}

\textit{Claim 3:} $\bar{\gamma}=\{(\gamma,\beta_1,\beta_2,...,\beta_a):\beta_i>\gamma\}$

The proof of Claim 3 is similar to the proof of \textit{Claim 3} in Theorem \ref{11322}.

\bigbreak

\label{1113222}

Then $\alpha_0$ is once again the least ordinal such that $G_{\alpha_0}=\bigcup_{\beta<\alpha_0}\bar{\beta}$.\qed

Similarly to Theorem \ref{11322}, we can identify the set of vertices beginning with each ordinal. We can also identify the set of vertices ending with each ordinal, as we can consider $N_v \cap G\backslash \bar{u}:u<v$ for every $v$ in $G$. However, we cannot in general which $(a+1)$-tuple ``generated'' some vertex $v$. For example, if $G(\omega,\sigma_{4,1})$, we cannot distinguish between say the vertices $(5,8,10,12,20)$ and $(5,9,11,12,20)$.

\bigbreak

\theorem Let  $\alpha$ and $\beta$ be infinite ordinals. Then \begin{align*}
G(\alpha,\sigma_{a,b})\cong G(\beta,\sigma_{a,b}) \Rightarrow \alpha = \beta
\end{align*}

\proof Omitted.

\bigbreak

\bigbreak

\remark Given the results above, one might  come to the conclusion that for any type $\tau$, $G(\alpha,\tau)\cong G(\beta,\tau) \Rightarrow \alpha = \beta$. However, it is easy to show that this is not the case. For example, $G(S,12)$ is simply the complete graph $K_{\kappa}$ where $\kappa$ is the cardinality of $S$. Thus, for example, $G(\omega+\omega,12)= G(\omega+1,12)=G(\omega,12)=K_\omega$.
It can be shown that for the type $123$, $G(\omega+\omega,123)= G(\omega+\omega+\alpha,123)$ for all $\alpha<\omega_1$ (this is left as an exercise to the reader).

It can also be shown that the graph $G(\alpha,3)$ is simply equal to a set of $\alpha$ isolated points, and that the graph $G(\alpha,\sigma_{0,b})$ is a set of ${\alpha\choose b}$ isolated points (also left as an exercise to the reader). This implies that we cannot necessarily distinguish between two graphs $G(\alpha,\sigma_{0,b})$ and $G(\alpha,\sigma_{0,b'})$, for example $G(\omega,3)\cong G(\omega,33)$.

\section{Chromatic Number}

We will now consider the chromatic number of $G(\alpha,\sigma_{a,b})$ in general. In \cite{shiftnew} it  is stated that Erd\"os and Hajnal showed that for any infinite cardinal $\kappa$, \begin{align*}
\chi(G(\kappa,132))=min\{\alpha:exp(\alpha)\geq \kappa\}
\end{align*} and that in general $\chi(G(\kappa,\sigma_{1,b}))=min\{\alpha:exp_{(b)}(\alpha)\geq \kappa\}$. This however covers just types of the form $\sigma_{1,b}=133...32$. Here we will give a fairly perspicuous proof that $\chi(G(\kappa,\sigma_{a,1}))=\kappa$ for $\kappa$ measurable, and a modification for $\kappa$ a strong limit. We will start with the measurable case, which is more direct.

\bigbreak

\definition Let $G$ be a graph. Then the \textit{chromatic number} of $G$, denoted by $\chi(G)$, is the minimal number of colours required to colour the vertices of $G$ so that no two neighbours  share the same colour.

\bigbreak

\remark Any clique of size $\kappa$ must be coloured with $\kappa$ many colours, and any co-clique of size $\kappa$ can be coloured with 1 colour. 

\bigbreak

\definition Let $\kappa$ be a cardinal. Then $\U$ is an \textit{ultrafilter} on $\kappa$ if:
\begin{enumerate}
\item $\U\subseteq \P(\kappa)$
\item $\emptyset \notin \U$
\item $x\cap y \in \U$ for $x,y \in \U$ (i.e. $\U$ is $\omega$-complete)
\item For each $x\in \P(\kappa)$ either $x\in \U$ or $\kappa \backslash x\in \U$
\end{enumerate}

\bigbreak

\definition $\U$ is \textit{principal} if there is some $\alpha<\kappa$ such that $\{\alpha\}\in \U$. 

This is a way for an ultrafilter to be trivial, as the sets in $\U$ are precisely the sets containing $\alpha$.

\bigbreak

\definition $\U$ is $\kappa$-complete if $\bigcap_{\alpha<\lambda} A_{\alpha} \in \U$ for $\lambda<\kappa$ and $A_{\alpha}\in \U$.

\bigbreak

Recall the following definition:

\definition An uncountable cardinal $\kappa$ is \textit{measurable} iff it has a $\kappa$-complete, non-principal ultrafilter.

\bigbreak

\theorem $\chi(G(\kappa,132))=\kappa$ for every measurable cardinal $\kappa$.

\proof Let $G=G(\kappa,132)$. Let $\bar{\alpha}$ be defined as the set of vertices beginning with $\alpha$ for each $\alpha<\kappa$.
Notice that given $\alpha$, for each $ \beta \geq \alpha$, $\bar{\alpha}$ contains a vertex which has a neighbour in $\bar{ \beta}$.  We would like to show that in our ``best-case scenario'', each $\bar{\alpha}$  is monochromatic with a distinct colour.

\center

\includegraphics[scale=0.3]{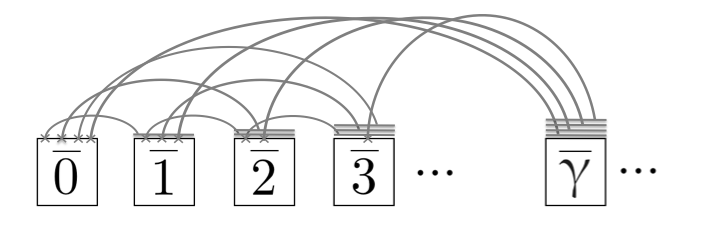}

\flushleft

 We will do this by showing that $G$ must have $\kappa$ many monochromatic sets, each contained in some $\bar{\alpha}$ and such that there is an edge from a vertex of each monochromatic set a vertex of every other monochromatic set, thus
resulting in a subgraph which behaves in a similar way to $K_{\kappa}$. We will construct such a subgraph by showing that there must be some ``large'' monochromatic subset of each $\bar{\alpha}$ (using the property that $\kappa$ is measurable).  Assume for a contradiction that $\chi(G)=\lambda<\kappa$.

Then $\bar{0}$ is coloured with $\lambda$ many colours, and so by the $\kappa$-completeness of the ultrafilter there is a monochromatic subset $S_0$ of $\bar{0}$ such that $T_0=\{t: (0,t)\in S_0\}\in \U$. Notice that $|T_0|=\kappa$. Let $t_0$ be the minimal element of $T_0$, and let $u_0=\{0\}$.

Now, $\bar{t_0}$ is coloured with $\lambda$ many colours, and so by the $\kappa$-completeness of the ultrafilter there is a monochromatic subset $S_1$ of $\bar{t_0}$ such that $\{t:(t_0,t)\in S_1\}\in \U$. Let $T_1=\{t:(t_0,t)\in S_1\}\cap T_0$, and notice that $T_1\in \U$. Thus $|T_1|=\kappa$ and $T_1\subseteq T_0$. Let $t_1$ be the minimal element of $T_1$, and let $u_1=\{0\} \cup \{t_0\}$.

Successor case: In general, suppose we have defined $T_{\alpha}$ for some $\alpha$ (and thus for all $\beta<\alpha$), and that $t_{\alpha}$ is the minimal element of $T_{\alpha}$. By the $\kappa$-completeness of the ultrafilter there is a monochromatic subset $S_{\alpha+1}$ of $\bar{t_{\alpha}}$ such that $\{t:(t_{\alpha},t)\in S_{\alpha+1}\}\in \U$. Let $T_{\alpha+1}=\{t:(t_{\alpha},t)\in S_{\alpha+1}\}\cap T_{\alpha}$, and notice that $T_{\alpha+1}\in \U$. Again $|T_{\alpha+1}|=\kappa$ and $T_{\alpha+1}\subseteq T_{\alpha}$. Let $t_{\alpha+1}$ be the minimal element of $T_{\alpha+1}$, and let $u_{\alpha+1}=\{0\}\cup\{t_\beta:\beta\leq \alpha\}$.

Limit case: For a limit ordinal $\gamma$, let $S_{\gamma}=\{(t_{\gamma},t):t\in T_{\gamma}\}$ where $T_{\gamma}=\bigcap_{\alpha<\gamma} T_{\alpha}$. Since $\U$ is $\kappa$-complete, $T_{\gamma}\in \U$. Let $t_{\gamma}$ be minimal in $T_{\gamma}$, and let $u_{\gamma}=\{0\}\cup\{t_\beta:\beta \leq \gamma\}$.

We repeat this process until $u=\cup_{\alpha<\kappa} u_\alpha$ has cardinality $\kappa$. This gives us $\kappa$ distinct monochromatic sets of vertices, each contained within some $\bar{\beta}$ for $\beta<\kappa$. Due to the construction of each $\bar{\beta}$, there is an edge between every pair of these monochromatic sets, and thus we will require $\kappa$ distinct colours, meaning   $\chi(G(\kappa,132))=\kappa$.\qed  

\label{chi132}

\bigbreak

\center

\includegraphics[scale=0.18]{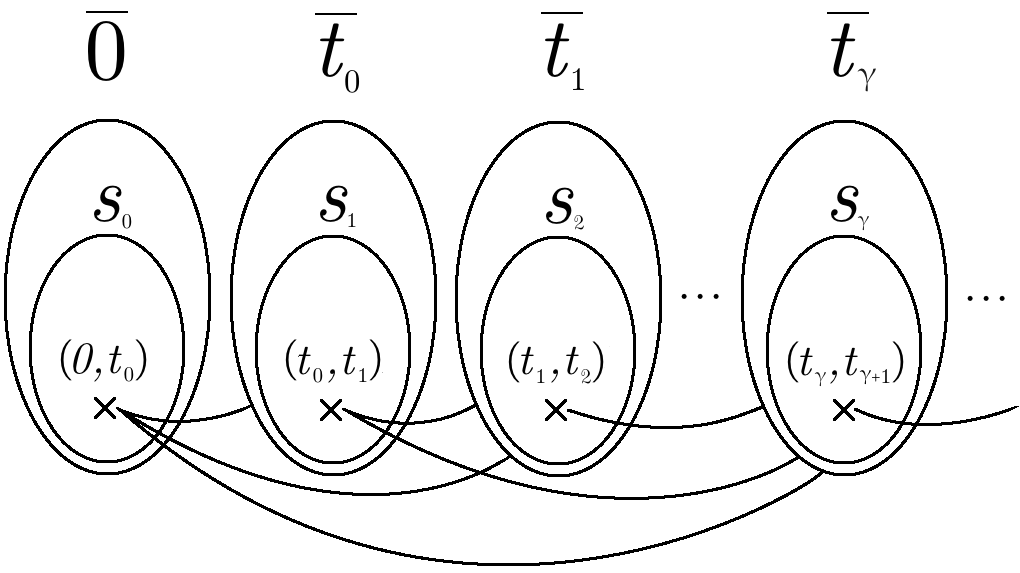}

\flushleft

\theorem $\chi(G(\kappa,\sigma_{a,1}))=\kappa$ for every measurable cardinal $\kappa$ and for every finite $a$.

\proof Let $G=G(\kappa,\sigma_{a,1})$. Let $\bar{\alpha}$ be defined as the set of vertices beginning with $\alpha$ for each $\alpha<\kappa$.  As in Theorem \ref{chi132}, we would like to show that in our ``best-case scenario'' each $\bar{\alpha}$ is monochromatic with a distinct colour. 

\flushleft

Notice that given $\alpha$, for each $ \beta \geq \alpha + a$, $\bar{\alpha}$ contains a vertex which has a neighbour in $\bar{ \beta}$. Once again, we will show that $G$ must have $\kappa$ many monochromatic sets, each contained in some $\bar{\alpha}$ and such that there is an edge from a vertex of each monochromatic set a vertex of every other monochromatic set, thus
resulting in a subgraph which behaves in a similar way to $K_{\kappa}$. We will construct such a subgraph by showing that there must be some ``large'' monochromatic subset of each $\bar{\alpha}$ (again, using the property that $\kappa$ is measurable). Assume for a contradiction that $\chi(G)=\lambda<\kappa$.

Then $\bar{0}$ is coloured with $\lambda$ many colours, and so by the $\kappa$-completeness of the ultrafilter there is a monochromatic subset $S_0$ of $\bar{0}$  such that $T_0=\{t: (0,1,2,...,a-1,t)\in S_0\}\in \U$. Notice that $|T_0|=\kappa$. Let $t_0$ be the minimal element of $T_0$, and let $u_0=\{0\}$.

The rest of the proof is similar to the proof of Theorem \ref{chi132}, except that at the successor stage we let $T_{\alpha+1}=\{t: (t_{\alpha},t_{\alpha}+1,t_{\alpha}+2,...,t_{\alpha}+a-1,t)\in S_{\alpha+1}\}\cap T_\alpha$,  and select $t_{\alpha+1}$ to be the minimal element of $T_{\alpha+1}$. Notice that $t_{\alpha+1}\geq t_{\alpha}+a$. \qed

\label{chi132b}

\bigbreak

Recall the following definition:
\definition A cardinal $\kappa$ is a \textit{strong limit} if  for all $\lambda<\kappa$, $2^\lambda<\kappa$.

\bigbreak

We can improve the above result by showing that it holds not just for measurable cardinals, but for all strong limit cardinals as well. We will require the following definition:
\definition Let $\lambda,\kappa$ be cardinals. Then $^\kappa[\lambda]$ is the set of all functions from $\kappa$ to $\lambda$.

\bigbreak

\theorem $\chi(G(\kappa,132))=\kappa$ for every strong limit cardinal $\kappa$.

\proof Assume for a contradiction that $\chi(G)=\lambda<\kappa$. In this proof, let every function mentioned  be strictly increasing.

 We will construct a tree $T$ where for each successor level the node $\bar{t}\in T$ (where $t\in \kappa$) is associated with both a vertex in the graph $G$ and some colour $< \lambda$. If $\bar{t}\in T$ lies on level $\delta$ and distinct neighbours $\bar{u}, \bar{v}$ of $\bar{t}$ lie on level $\delta+1$, we have that $(t,u)$ in $G$ has some colour $c_1$ and $(t,v)$ in $G$ has some colour $c_2\neq c_1$. Additionally, $u$ is the minimal point in $T_1\subseteq \kappa$ such that $(t,u)$ has colour $c_1$ and $v$ is the minimal point in $T_2\subseteq \kappa$ such that $(t,v)$ has colour $c_2$.

The root of $T$ is $\bar{0}$, which is not associated with any particular colour. We will now describe the first level of $T$. Since $\bar{0}$  is coloured with at most $\lambda$ many colours, $\bar{0}$ is joined to a node on level 1 corresponding to each colour which arises. Moreover, we would like each node to be the smallest element of $\kappa$ corresponding to each colour. First, we partition $\kappa$ into disjoint sets $\{T_{\alpha}:\alpha<\lambda\}$ (some of which may be empty), where $(0,t)$ is coloured with colour $c_{\alpha}$ for all $t\in T_{\alpha}$. We then pick the smallest non-empty element of each $T_{\alpha}$, and call it $t_{\alpha}$. Level 1 of $T$ is thus  $\{\bar{t_0}, \bar{t_1}, \bar{t_2}, \bar{t_3}, ...,\bar{t_\alpha},...\}$ where $\alpha<\lambda$. We can write each element of level 1 of $T$ as $\bar{t_{f}}$ where $f$ is a function from $1$ to $\lambda$, i.e. $f\in {^{1}[\lambda]}$, and let the corresponding subset of $\kappa$ be $T_f$.

We will define each node in a successor level of $T$ in a similar manner. Given $\bar{t_{f}}$ on level $\delta$ with $f\in{^\delta[\lambda]}$, when considering $\bar{t_{g}}$  with $g\in{^{\delta+1}[\lambda]}$ on level $\delta+1$ extending $f$ we partition $T_{f}$ into disjoint sets $T_{g}$ where $(t_{f},t)$ is coloured with colour $c_{g}$ for each $t\in T_{g}$. Notice that $T_{g}\subseteq T_{f}$. 

For each $f\in {^\gamma[\lambda]}$ and limit $\gamma$, we will now define the node $\bar{t_f}$ on level $\gamma$ of $T$.  We define $T_f$ to be the intersection of the decreasing sequence $(T_{f|\xi}: \xi < \gamma)$, and if $T_f \neq \emptyset$ we choose $t_f$ to be its least member. By construction, $T$ must have size $\kappa$, and so such a $t_f$ exists for each limit level $\gamma<\kappa$. Since $\kappa$ is a strong limit, the width of $T$ at each level $\delta<\kappa$ will be at most $\lambda^\delta<\kappa$, and thus $T$ must have height $\kappa$.

\center

\includegraphics[scale=0.2]{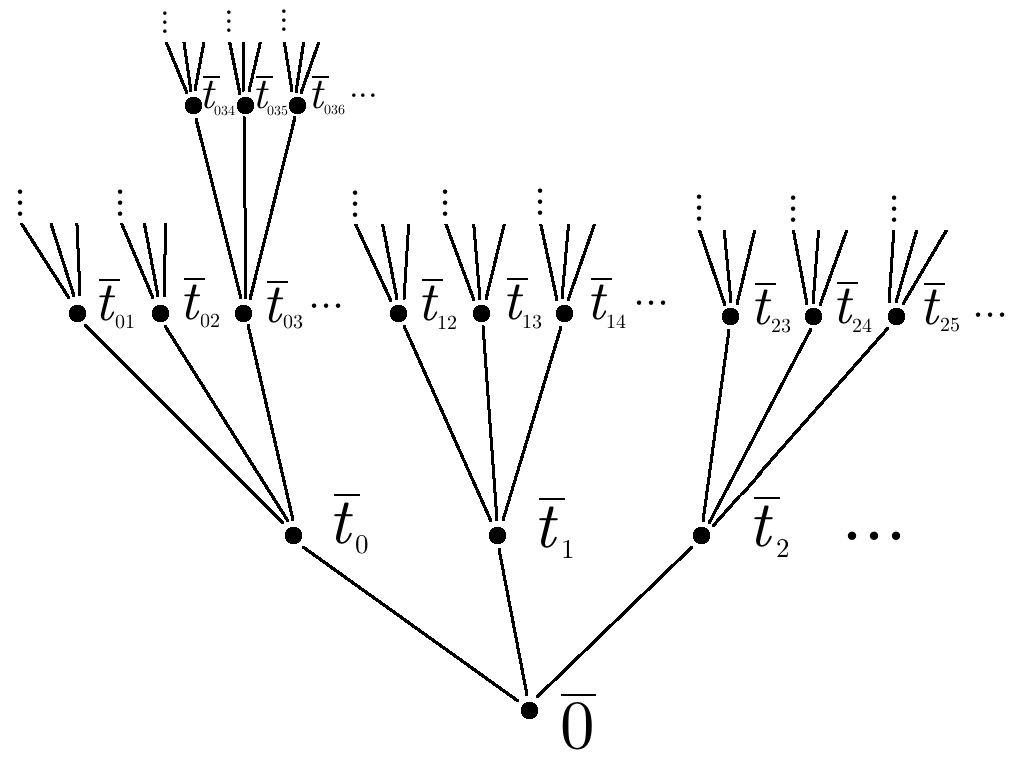}

\flushleft

Each branch in $T$ is similar to the construction in Theorem \ref{chi132}. As $T$ has height $\kappa>\lambda$ (and $\kappa$ is a strong limit), and the height of a tree is the supremum of the length of all its branches, there must be some branch in $T$ of length $\lambda^+ >\lambda$. Due to the construction of each $\bar{0}, \bar{1}, ..., \bar{\beta}$ and of our tree $T$, we have $\lambda^+$ distinct monochromatic sets of vertices with an edge between every pair of these sets of vertices, and thus $G$ is coloured with $\lambda^+>\lambda$ distinct colours, a contradiction. Thus $\chi(G(\kappa,132))=\kappa$.\qed

\label{compact}

\bigbreak

\bigbreak

\theorem $\chi(G(\kappa,\sigma_{a,1}))=\kappa$ for every strong limit cardinal $\kappa$ and every finite $a$.

\proof Assume for a contradiction that $\chi(G)=\lambda<\kappa$. We will show that $\chi(G)>\lambda$. In this proof, let every function mentioned be strictly increasing.

We will construct a  tree $T$ similar to the one in Theorem $\ref{compact}$. The root is $\bar{0}$, which is not associated with any particular colour. Now, since $\bar{0}$ is coloured with at most $\lambda$ many colours, we let each node on level 1 of $T$ correspond to the minimal vertex $t_{\alpha}$ in $\kappa$ such that $(0,1,...,a-1,t_{\alpha})$ is coloured with colour $c_{\alpha}$. We  denote the vertices in level 1 of $T$ by $\bar{t_{f}}$ where $f\in {^1[\lambda]}$, and let the corresponding subset of $\kappa$ be $T_f$.
In general, given $\bar{t_{f}}$ on level $\delta$ with $f\in{^\delta[\lambda]}$, when considering $\bar{t_{g}}$  with $g\in{^{\delta+1}[\lambda]}$ on level $\delta+1$ extending $f$ we partition $T_{f}$ into disjoint sets $T_{g}$ where $(t_{f},t_f+1,...,t_f+a-1,t)$ is coloured with colour $c_{g}$ for each $t\in T_{g}$. Notice that $T_{g}\subseteq T_{f}$.

For each $f\in {^\gamma[\lambda]}$ and limit $\gamma$, we will now define the node $\bar{t_f}$ on level $\gamma$ of $T$.  We define $T_f$ to be the intersection of the decreasing sequence $(T_{f|\xi}: \xi < \gamma)$, and if $T_f \neq \emptyset$ we choose $t_f$ to be its least member. If $t_f$ with $f\in{^\gamma[\lambda]}$ is the minimal point in this intersection, we let the corresponding node at level $\gamma$ of $T$ be $\bar{t_f}$, and continue as before. By the same argument as in Theorem \ref{compact}, $T$ must have height $\kappa$ and width $<\kappa$ at each level. Thus there must be some branch of $T$ of length $\lambda^+>\lambda$, and so $G$ is coloured with $\lambda^+$ distinct colours, a contradiction. Thus $\chi(G(\kappa,\sigma_{a,1}))=\kappa$.\qed

\bigbreak

\pagebreak

\section{Automorphism Groups}

In addition to recognising the ordinal $\alpha$ inside $G(\alpha,\tau)$, we are also interested in the problem of recognising the graph $G$ from its automorphism group. This problem has been considered by M. Rubin for many classes of structures, for instance see \cite{rubin}. 

We will see that in general, the graph $G(\alpha,\sigma_{a,b})$ cannot be constructed from its automorphism group (see Propositions \ref{aut132inf} and \ref{aut132fin}). However, we can construct $G(n,\sigma_{a,1})$ from its automorphism group for finite $n$ (where $n$ is sufficiently large) and $a\geq 3$ (see Theorem \ref{autgroupb1}).

\bigbreak

\proposition The automorphism group of $G(\alpha,132)$ for any infinite ordinal $\alpha$ is the trivial group $\{id\}$.

\proof By the proof of Theorem \ref{132iso} \textit{Claim 3},  we see that $G(\alpha,132)$ has only the trivial  automorphism.\qed

\label{aut132inf}

\bigbreak

\proposition The automorphism group of $G(n,132)$ for finite $n$ is $\Z_2$.

\proof Using the same method as in Theorem \ref{132iso}, we can identify each element in $G(n,132)$ up to order reversal of $n$. Thus the automorphism group of $G(n,132)$ is $\Z_2$.\qed

\label{aut132fin}

\bigbreak

Similarly, by Theorem \ref{isomorphism}, we see that the automorphism group of $G(n,\sigma_{1,b})$ for  $n,b<\omega$ is $\Z_2$ and the automorphism group of $G(\alpha,\sigma_{1,b})$ for any infinite ordinal $\alpha$ is the trivial group $\{id\}$.

\bigbreak

\proposition Call two vertices \textit{equivalent} if they are not neighbours and they have the exact same neighbour set, denoted by $u\sim v$ (this is an equivalence relation). Then if two vertices in a graph are equivalent, there exists an automorphism permuting exactly these two points (and fixing everything else).

\proof The function which permutes $u,v$ and fixes everything else is an isomorphism.\qed

\label{equivalence}

\bigbreak

\theorem The automorphism group of $G(n,11322)$ for finite $n\geq 5,\ n\neq 7$ is:
\begin{align*}
\Z_2\times S_{4(n-3)} \times \prod_{j=1}^{n-5} (S_{2j+1})^2  \times \prod_{j=1}^{n-6} (S_{n-5-j})^j
\end{align*}

\proof 
We see by definition of $G(n,11322)$ that two vertices $x=(x_1,x_2,x_3)$ and $y=(y_1,y_2,y_3)$  are equivalent iff \begin{align*}
[\text{either } x_1=y_1 \text{ or } x_1,y_1\in \{0,1\}] \text{ and }
[\text{either } x_3=y_3 \text{ or } x_3,y_3\in \{n-1,n-2\}]
\end{align*}
Moreover, these equivalent vertices can be permuted independently of one another by Proposition \ref{equivalence}. The automorphism group therefore contains  the direct product of a set of permutation groups of these equivalence classes, which we will split into the following 4 cases:
\begin{enumerate}[(1)]
\item  Isolated points, i.e. vertices $(x, y, z)$ with $x\in \{0, 1\}$ and $z\in \{n-1, n-2\}$
\item Vertices $(x, y, z)$ with $x\in \{0,1\}$ and $z\notin \{n-1, n-2\}$
\item Vertices $(x, y, z)$ with $x\notin \{0,1\}$ and $z\in \{n-1, n-2\}$
\item   Vertices $(x, y, z)$ with $x\notin \{0,1\}$ and $z\notin \{n-1, n-2\}$
\end{enumerate}

We will then collapse each equivalence class to a point and show that for $n\geq 8$ the resulting graph has exactly one non-trivial automorphism, namely the order reversal automorphism mapping each point $(x_1, x_2, x_3)$ to $(n-1-x_3, n-1-x_2, n-1-x_1)$. 

\bigbreak

(1) \underline{{Isolated Points}}

The ordinal $n=\{0,1,2,...,n-1\}$. The isolated points are the vertices starting with $0$ or $1$ and ending with $n-1$ or $n-2$. There are:

$\bullet$ $n-2$ vertices starting with 0 and ending with $n-1$,

$\bullet$ $n-3$ vertices starting with 1 and ending with $n-1$,

$\bullet$ $n-3$ points starting with 0 and ending with $n-2$, and 

$\bullet$ $n-4$ points starting with 1 and ending with $n-2$ ,

giving a total of $4(n-3)$.  Thus we have $4(n-3)$ vertices that can be permuted freely, meaning that  $S_{4(n-3)}$ lies in the direct product.\newline

(2) \underline{Vertices $(x,y,z)$ with $x\in \{0,1\}$ and $z\notin \{n-1, n-2\}$, or}

\center
(3) \underline{Vertices $(x,y,z)$ with $x\notin \{0,1\}$ and $z\in \{n-1,n-2\}$}
\flushleft

Now consider Case (2), i.e. the vertices starting with either $0$ or $1$, and ending with $n-j-1$ for $j>1$. There are $n-j-2$ vertices beginning with 0 and ending with $n-j-1$, and $n-j-3$ vertices beginning with 1 and ending with $n-j-1$, giving a total of $2n-2j-5$. Now, $j\geq 2$, meaning we have $2n-9$ vertices beginning with 0 or 1 and ending in $n-3$, $2n-11$ vertices beginning with 0 or 1 and ending in $n-4$, $2n-13$ vertices beginning with 0 or 1 and ending in $n-5$, and so on. Thus the Automorphism group ``includes'' $S_{2n-9} \times S_{2n-11}\times S_{2n-13} \times ...$. We need to have two copies of each of these as we need to account for Case (3) which is similar to Case (2), and so $\prod_{j=1}^{n-5} (S_{2j+1})^2$ lies in the direct product. \newline

(4) \underline{{Vertices $(x,y,z)$ with $x\notin \{0,1\}$ and $z\notin \{n-1,n-2\}$}}

Finally, we consider the vertices beginning with $x>1$ and ending with $z<n-2$. First consider the vertices beginning with 2. 

$\bullet$ There are $n-7$ vertices beginning with 2 and ending with $n-3$,

$\bullet$ There are $n-8$ vertices beginning with 2 and ending with $n-4$,

$\bullet$ There are $n-9$ vertices beginning with 2 and ending with $n-5$, and so on.

Similarly,

$\bullet$ There are $n-8$ vertices beginning with 3 and ending with $n-3$,

$\bullet$ There are $n-9$ vertices beginning with 3 and ending with $n-4$,

$\bullet$ There are $n-10$ vertices beginning with 3 and ending with $n-5$, and so on.

In general, there are $n-(j+k+3)$ vertices beginning with $j$ and ending with $n-k-1$.

Thus  $S_{n-7}\times (S_{n-8})^2\times  (S_{n-9})^3  \times ... \times {(S_{n-k})}^{(k-6)}$  i.e.
$\prod_{j=1}^{n-6} (S_{n-5-j})^j$ lies in the direct product.\newline

\underline{{Equivalence Graph}}

We will now collapse all the equivalence classes in $G(n,11322)$ to form a new graph $G'$, and show that the only non-trivial automorphism $G'$ has is the order reversal automorphism (except for some special cases where $n<8$).

Now, to collapse the equivalence classes we adopt the following notation:

$\bullet$ Isolated vertices $(x,y,z)$ with $x\in \{0,1\}$ and $z\in \{n-1,n-2\}$ are denoted by $[1,n-2]$.

$\bullet$  Vertices $(x,y,z)$ with $x\in \{0,1\}$ and $z\notin \{n-1,n-2\}$ are denoted by $[1,z]$

$\bullet$ Vertices $(x,y,z)$ with $x\notin \{0,1\}$ and $z\in \{n-1,n-2\}$ are denoted by $[x,n-2]$

$\bullet$   Vertices $(x,y,z)$ with $x\notin \{0,1\}$ and $z\notin \{n-1,n-2\}$ are denoted by $[x,z]$

$\bullet$ The vertices $(0,1,2)$ and $(n-3,n-2,n-1)$ are taken as special cases and denoted by $[0,2]$ and $[n-3,n-1]$ respectively.

Notice that $G'$ is a subgraph of the graph $G(n,132)$. $G'$ has exactly one isolated point, namely $[1,n-2]$.

\bigbreak

We now define $\bar{\xi_0}$ to be the following: $\bar{\xi_0}$ is a set of vertices in $G'$ such that:

$\circ$ The induced subgraph on $\bar{\xi_0}$ is a maximal co-clique

$\circ$ All neighbour sets of vertices in the induced subgraph on  $\bar{\xi_0}$ are pairwise disjoint

\textit{Claim 1:} The only sets of vertices in $G'$ satisfying the conditions for $\bar{\xi_0}$ are $\{[1,x]:x>1\}\cup\{[0,2]\}$ and $\{[y,n-2]:y<n-2\}\cup\{[n-3,n-1]\}$.

\textit{Proof of Claim 1:} First we see that both $\{[1,x]:x>1\}\cup\{[0,2]\}$ and $\{[y,n-2]:y<n-2\}\cup\{[n-3,n-1]\}$ satisfy the conditions above. We now need to show that these are the \textit{only} sets satisfying the conditions above.

Suppose $\bar{\xi_0}\not\subseteq \{[1,x]:x>1\}\cup\{[0,2]\}$. Then $\bar{\xi_0}$ contains some vertex $[a,b]$ with $a\notin \{0,1\}$. Then $[0,b]$ must lie in  $\bar{\xi_0}$ as $\bar{\xi_0}$ is maximal; but then either the neighbour sets of $\bar{\xi_0}$ are not pairwise disjoint, or $b\in \{n-2,n-1\}$ in which case $[a,b]\in \{[y,n-2]:y<n-2\}\cup\{[n-3,n-1]\}$. 

We now see that every other vertex in $\bar{\xi_0}$ must lie in  $\{[y,n-2]:y<n-2\}\cup\{[n-3,n-1]\}$. If not, there is some vertex $[c,d]$ in $\bar{\xi_0}$ with $d\notin \{n-2,n-1\}$. If $d<a$ then $[d,a]$ is a neighbour of both $[a,b]$ and $[c,d]$, contradicting the second condition above. If $d=a$ we no longer have a co-clique. Thus we must have $d>a$; but then $[d,a]$ must also lie in $\bar{\xi_0}$ as $\bar{\xi_0}$ is maximal, and thus the neighbour sets of $\bar{\xi_0}$  are not pairwise disjoint (namely, $[c,a]$ and $[d,a]$ share at least one neighbour, and so do $[d,a]$ and $[d,b]$). Hence $\bar{\xi_0}= [a,b]\in \{[y,n-2]:y<n-2\}\cup\{[n-3,n-1]\}$

Similarly, we see that if $\bar{\xi_0}\not\subseteq [a,b]\in \{[y,n-2]:y<n-2\}\cup\{[n-3,n-1]\}$, then $\bar{\xi_0}= \{[1,x]:x>1\}\cup\{[0,2]\}$. 

Thus either $\bar{\xi_0}=\{[1,x]:x>1\}\cup\{[0,2]\}$ or $\bar{\xi_0}=\{[y,n-2]:y<n-2\}\cup\{[n-3,n-1]\}$.
\flushright
 \textit{Q.E.D. \ \ Claim 1} 
\flushleft

\bigbreak

We now define $\bar{\xi_1}$ to be the following: $\bar{\xi_1}$ is a set of vertices in $G'\backslash \bar{\xi_0}$ such that:

$\circ$ The induced subgraph on $\bar{\xi_1}$ is a maximal co-clique

$\circ$ All neighbour sets of vertices in the induced subgraph on  $\bar{\xi_1}$ are pairwise disjoint

$\circ$ Exactly one vertex in the induced subgraph on $\bar{\xi_0}$ has neighbour set $\bar{\xi_1}$

\textit{Claim 2:}  If $\bar{\xi_0}=\{[1,x]:x>1\}\cup\{[0,2]\}$ then $\bar{\xi_1}=\{[2,x]:x>2\}$, and if $\bar{\xi_0}=\{[y,n-2]:y<n-2\}\cup\{[n-3,n-1]\}$ then $\bar{\xi_1}=\{[y,n-3]:y<n-3\}$.

\textit{Proof of Claim 2:} Suppose $\bar{\xi_0}=\{[1,x]:x>1\}\cup\{[0,2]\}$. If exactly one vertex in $\bar{\xi_0}$ has neighbour set $\bar{\xi_1}$, we must have that $\bar{\xi_1}=\{[z,x]:x>z\}$ for some $z>1$. If $z>2$, then the vertex $[2,b]$ must also lie in $\bar{\xi_1}$ for some $z<b<n-2$, a contradiction as the neighbour sets would not be pairwise disjoint (namely, the vertex $[b,b+2]$ is neighbours with both $[2,b]$ and $[z,b]$).

Now suppose $\bar{\xi_0}=\{[y,n-2]:y<n-2\}\cup\{[n-3,n-1]\}$. Again, If exactly one vertex in $\bar{\xi_0}$ has neighbour set $\bar{\xi_1}$, we must have that $\bar{\xi_1}=\{[y,z]:y<z\}$ for some $z<n-2$. If $z<n-3$, then the vertex $[a,n-3]$ must also lie in $\bar{\xi_1}$ for some $2<a<z$, a contradiction as the neighbour sets would not be pairwise disjoint (namely, the vertex $[a-2,a]$ is neighbours with both $[a,n-3]$ and $[a,z]$).
\flushright
 \textit{Q.E.D. \ \ Claim 2} 
\flushleft

In general, for $i\geq 2$ we define $\bar{\xi_i}$ to be the following: $\bar{\xi_i}$ is a set of vertices in $G'\backslash \bigcup_{j<i}\bar{\xi_j}$ such that:

$\circ$ The induced subgraph on $\bar{\xi_i}$ is a maximal co-clique

$\circ$ All neighbour sets of vertices in the induced subgraph on  $\bar{\xi_i}$ are pairwise disjoint

$\circ$ Exactly one vertex in the induced subgraph on $\bar{\xi_{i-2}}$ has neighbour set $\bar{\xi_i}$

\textit{Claim 3:}  If $\bar{\xi_0}=\{[1,x]:x>1\}\cup\{[0,2]\}$ then $\bar{\xi_i}=\{[i-1,x]:x>i-1\}$, and if $\bar{\xi_0}=\{[y,n-2]:y<n-2\}\cup\{[n-3,n-1]\}$ then $\bar{\xi_i}=\{[y,n-i-2]:y<n-i-2\}$.

\textit{Proof of Claim 3:} Suppose $\bar{\xi_0}=\{[1,x]:x>1\}\cup\{[0,2]\}$. If exactly one vertex in $\bar{\xi_{i-2}}$ has neighbour set $\bar{\xi_i}$, we must have that $\bar{\xi_i}=\{[z,x]:x>z\}$ for some $z>i-2$. If $z>i-1$, then the vertex $[i-1,b]$ must also lie in $\bar{\xi_i}$ for some $z<b<n-2$, a contradiction as the neighbour sets would not be pairwise disjoint (namely, the vertex $[b,b+2]$ is a neighbour of  both $[i-1,b]$ and $[z,b]$).

Now suppose $\bar{\xi_0}=\{[y,n-2]:y<n-2\}\cup\{[n-3,n-1]\}$. Again, If exactly one vertex in $\bar{\xi_{i-2}}$ has neighbour set $\bar{\xi_i}$, we must have that $\bar{\xi_i}=\{[y,z]:y<z\}$ for some $z<n-i-1$. If $z<n-i-2$, then the vertex $[a,n-i-2]$ must also lie in $\bar{\xi_i}$ for some $2<a<z$, a contradiction as the neighbour sets would not be pairwise disjoint (namely, the vertex $[a-2,a]$ is a neighbour of  both $[a,n-3]$ and $[a,z]$).
\flushright
 \textit{Q.E.D. \ \ Claim 3} 
\flushleft

Thus the automorphism group for $G(n,11322)$ where $n\geq 5,\ n\neq 7$ is
\begin{align*}
\Z_2\times S_{4(n-3)} \times \prod_{j=1}^{n-5} (S_{2j+1})^2 \times \prod_{j=1}^{n-6} (S_{n-5-j})^j
\end{align*}\qed

\label{finiteaut11322}

In the case where $n=7$, the graph has an extra connected component, and thus the automorphism group is 
\begin{align*}
\Z_2\times \Z_2\times S_{4(n-3)} \times \prod_{j=1}^{n-5} (S_{2j+1})^2 \times \prod_{j=1}^{n-6} (S_{n-5-j})^j
\end{align*}

\flushleft

However, for $n\geq 5,\ n\neq 7$, there is only one connected component (apart from the isolated point). Hence there is only one non-trivial automorphism of $G'$, namely the order reversing automorphism which means that $\Z_2$ lies in the direct product which makes up the automorphism group of $G(n,11322)$.

\bigbreak

\bigbreak

We can rewrite the automorphism group in Theorem \ref{finiteaut11322} as follows:
\begin{align*}
\Z_2\times S_{\big [{n-2 \choose 1} + 2\cdot {n-3 \choose 1} + \cdot {n-4 \choose 1}  \big ]} \times \prod_{j=1}^{n-5} \left(S_{{n-j-1\choose 1} + {n-j-2 \choose 1}}\right)^2 \times \prod_{j=1}^{n-6} \left(S_{{n-5-j\choose 1}}\right)^j
\end{align*}

It may seem like we are overcomplicating things here, but writing the automorphism group like this makes it easier to see how this case relates to the general case.

The following results  are given without proof, and are based on having a good enough understanding of the graph structure to see which parts can or cannot be permuted ``independently''.

\bigbreak

\theorem The automorphism group of $G(n,1113222)$ for finite $n\geq 12$ is:
\begin{align*}
\Z_2\times S_{\big [{n-2 \choose 2} + 2\cdot {n-3 \choose 2} + 3\cdot {n-4 \choose 2} + 2\cdot {n-5 \choose 2} +{n-6 \choose 2} \big ] } \times \prod_{j=4}^{n-4} \left(S_{\big [{n-j-1 \choose 2} + {n-j-2 \choose 2} + {n-j-3 \choose 2}\big ]}\right)^2  \times \prod_{j=1}^{n-8} \left(S_{{n-7-j \choose 3}}\right)^j
\end{align*}

\label{finiteaut1113222}

\bigbreak

\bigbreak

\theorem The automorphism group of $G(n, \sigma_{a,1})$ for $a\geq 3$ and $n$ sufficiently large is:
\begin{align*}
\Z_2 &\times S_{\sum_{k=1}^{a-1} k\cdot [{n-k-1 \choose {a-1}} + {n-2a+k-1 \choose {a-1}}]+ a\cdot {n-a-1\choose {a-1}} }\\ &\times \prod_{j={a+1}}^{n-a-1} \left(S_{\sum_{k=0}^{a-1} {n-j-k-1 \choose {a-1}}}\right)^2  \times \prod_{j=1}^{n-2(a+1)} \left(S_{{n-2a-j-1 \choose 3}}\right)^j
\end{align*}

\label{autgroupb1}

\bigbreak

\bigbreak

We can now consider some infinite cases.

\bigbreak

\theorem The automorphism group of $G(\omega,11322)$ is $\prod_{n=1}^{\omega}(S_n)^{\omega}$

\textit{Sketch of Proof.} The points which cannot be identified are the points $(x,y,z)$ where $x$ and $z$ are fixed, and thus the number of values $y$ can take defines how many permutations we can perform on this subgroup.

(For example, if $x=3$ and $z=8$ then $y$ can take 4 values, illustrating how the points $(3,4,8), (3,5,8), (3,6,8), (3,7,8)$ can all be permuted.)

Now, $y$ can take exactly $n<\omega$ values (for each $n$) in $\omega$ many ways, and so the automorphism group of $G(\omega,11322)$ is $\prod_{n=1}^{\omega}(S_n)^{\omega}$.\qed

\bigbreak

\theorem The automorphism group of $G(\omega+1,11322)$ is 
\begin{align*}
\prod_{n=1}^{\omega}(S_n)^{\omega} \times (S_\omega)^{\omega}
\end{align*}

\textit{Sketch of Proof.} The points which cannot be identified are the points $(x,y,z)$ where $x$ and $z$ are fixed, and thus the number of values $y$ can take defines how many permutations we can perform on this subgroup. 

Now, $y$ can take exactly $n$ values (for each $n<\omega$) in $\omega$ many ways, and additionally $y$ can take $\omega$ many values in $\omega$ many ways (consider vertices of the form $(x,y,\omega)$) and so the automorphism group of $G(\omega,11322)$ is $\prod_{n=1}^{\omega}(S_n)^{\omega}\times (S_\omega)^{\omega}$.\qed

\bigbreak

\theorem For every countable ordinal $\alpha>\omega$, the automorpnism group of $G(\alpha,11322)$ is: 
\begin{align*}
\prod_{n=1}^{\omega}(S_n)^{\omega} \times (S_\omega)^{\omega}
\end{align*}

\textit{Sketch of Proof.} 
The points which cannot be identified are the points $(x,y,z)$ where $x$ and $z$ are fixed, and thus the number of values $y$ can take defines how many permutations we can perform on this subgroup. 

Now, $y$ can take exactly $n$ values (for each $n<\omega$) in $\omega$ many ways, but $y$ can still only take $\omega$ many values in $\omega$ many ways (consider vertices of the form $(x,y,\alpha)$ where $\alpha=\alpha'+k$ for some limit $\alpha'\geq \omega$ and $x<\alpha'$) and so the automorphism group of $G(\alpha,11322)$ is $\prod_{n=1}^{\omega}(S_n)^{\omega}\times (S_\omega)^{\omega}$.\qed

\bigbreak

\theorem The automorphism group of $G(\omega,1113222)$ is $\prod_{n=2}^{\omega}\left(S_{n \choose 2}\right)^{\omega}$

\textit{Sketch of Proof.} The points which cannot be identified are the points $(x,y_1,y_2,z)$ where $x$ and $z$ are fixed, and thus the number of values $y_1$ and $y_2$ can take defines how many permutations we can perform on this subgroup. 

Now, $y_1$ and $y_2$ can take exactly ${n\choose 2}$ values (for each $n<\omega$) in $\omega$ many ways, and so the automorphism group of $G(\omega,1113222)$ is  $\prod_{n=2}^{\omega}\left(S_{n \choose 2}\right)^{\omega}$.\qed

\bigbreak

We can generalise this infinite case to the following:

\theorem The automorphism group of $G(\omega, \sigma_{a,1})$ is $\prod_{n=(a-1)}^{\omega}\left(S_{n \choose a-1}\right)^{\omega}$

\textit{Sketch of Proof.} The points which cannot be identified are the points $(x,y_1,y_2,...,y_{a-1},z)$ where $x$ and $z$ are fixed, and thus the number of values the $y_i$'s  can take defines how many permutations we can perform on this subgroup. 

Now, the $y_i$'s can take exactly ${n\choose {a-1}}$ values (for each $n<\omega$) in $\omega$ many ways, and so the automorphism group of $G(\omega, \sigma_{a,1})$ is  $\prod_{n=(a-1)}^{\omega}\left(S_{n \choose a-1}\right)^{\omega}$.\qed

\bigbreak

\theorem For any countable $\alpha>\omega$, the automorphism group of $G(\alpha, \sigma_{a,1})$   is $\prod_{n=(a-1)}^{\omega}\left(S_{n \choose a-1}\right)^{\omega}\times (S_\omega)^{\omega}$.

\textit{Sketch of Proof.} The points which cannot be identified are the points $(x,y_1,y_2,...,y_{a-1},z)$ where $x$ and $z$ are fixed, and thus the number of values the $y_i$'s  can take defines how many permutations we can perform on this subgroup. 

Now, the $y_i$'s can take exactly ${n\choose {a-1}}$ values (for each $n<\omega$) in $\omega$ many ways, and $\omega$ many values in $\omega$ many ways (consider vertices of the form $(x,y_1,y_2,...,y_{a-1},\alpha)$ where $\alpha=\alpha'+k$ for some limit $\alpha'\geq \omega$ and $x<\alpha'$). and so the automorphism group of $G(\alpha, \sigma_{a,1})$ is  $\prod_{n=(a-1)}^{\omega}\left(S_{n \choose a-1}\right)^{\omega}\times (S_\omega)^{\omega}$.\qed

\bigbreak

\theorem For any cardinal $\kappa$, the automorphism group of $G(\kappa, \sigma_{a,b})$  is $\prod_{n=(a-1)}^{\omega}\left(S_{n \choose a-1}\right)^{\omega}\times \prod_{\gamma < \kappa}(S_\gamma)^{\gamma}$.

\bigbreak

\theorem For any ordinal $\alpha$ of cardinality $\kappa$ with $\alpha>\kappa$, the automorphism group of $G(\alpha, \sigma_{a,b})$  is $\prod_{n=(a-1)}^{\omega}\left(S_{n \choose a-1}\right)^{\omega}\times \prod_{\gamma\leq \kappa}(S_\gamma)^{\gamma}$.

\bigbreak

\section*{Final Remarks}

Our goal is to extend these results to all $G(S,\tau)$ for every totally ordered set $S$ and type $\tau$, and also to their automorphism groups. We would also like to find the chromatic number of these graphs in general, and to extend the tree construction used in Section 4 to consider colourings of more shift graphs.

\pagebreak

\end{document}